\documentclass[review]{elsarticle}
\usepackage{lipsum}
\makeatletter
\def\ps@pprintTitle{%
 \let\@oddhead\@empty
 \let\@evenhead\@empty
 \def\@oddfoot{}%
 \let\@evenfoot\@oddfoot}
 
\usepackage{lineno,hyperref}
\modulolinenumbers[5]

\journal{Journal of Computational Science}

\usepackage{color}
\usepackage{graphicx}
\usepackage{epstopdf}
\usepackage{amsmath}
\usepackage{amssymb}
\usepackage{subfigure}
\usepackage{float}

\newcommand{\bv}[1]{\boldsymbol{#1}}
\newcommand{\domain}{\Omega}
\newcommand{\bound}{\Gamma}
\newcommand{\R}{\mathbb{R}}
\newcommand{\Cnn}{\mathbb{C}^{N\times N}}
\newcommand{\Cn}{\mathbb{C}^N}
\newcommand{\bb}{\bv{b}}

\newcommand{\sect}{Section }

\newcommand{\RNum}[1]{\uppercase\expandafter{\romannumeral #1\relax}}

\begin{document}

\begin{frontmatter}

\title{A Matrix-free Preconditioner for the Helmholtz Equation based on the Fast Multipole Method}

\author[kaust]{Huda Ibeid\corref{corresponding}}
\author[tokyo]{Rio Yokota}
\author[kaust]{David Keyes}


\address[kaust]{Extreme Computing Research Center, King Abdullah University of Science and Technology, Thuwal, Saudi Arabia}
\address[tokyo]{Global Scientific Information and Computing Center, Tokyo Institute of Technology, Tokyo, Japan}

\begin{abstract}
Fast multipole methods (FMM) were originally developed for accelerating $N$-body problems for particle-based methods. FMM is more than an $N$-body solver, however. Recent efforts to view the FMM as an elliptic Partial Differential Equation (PDE) solver have opened the possibility to use it as a preconditioner for a broader range of applications. FMM can solve Helmholtz problems with optimal $\mathcal{O}(N \log N)$ complexity, has compute-bound inner kernels, and highly asynchronous communication patterns. The combination of these features makes FMM an interesting candidate as a preconditioner for sparse solvers on architectures of the future. The use of FMM as a preconditioner allows us to use lower order multipole expansions than would be required as a solver because individual solves need not be accurate. This reduces the amount of computation and communication significantly and makes the time-to-solution competitive with state-of-the-art preconditioners. Furthermore, the high asynchronicity of FMM allows it to scale to much larger core counts than factorization-based and multilevel methods. We describe our tests in reproducible details with freely available codes.
\end{abstract}


\end{frontmatter}


\section{ Introduction}

The Helmholtz equation can be used to describe both wave propagations and scattering phenomena arising in many fields of science and technology. For example, Helmholtz equations are used to express acoustic phenomena in aeronautics~\cite{Laird2002} and underwater acoustics~\cite{Arnold1998,Fix1978}. They are also utilized in geophysical applications~\cite{Plessix2004} and in electromagnetic applications, e.g., photolithography~\cite{Urbach1991}. Nevertheless, the development of accurate, robust, and efficient numerical methods for the solution of the Helmholtz equation with high wavenumber, and therefore highly oscillatory solution, remains an important challenge~\cite{zienkiewicz2000}.

One important aspect in the numerical methods for solving partial differential equations (PDEs) is the efficient solution of the large, and usually sparse, linear systems arising from the PDE discretization. In particular,  we look for solutions of the Helmholtz equation discretized using finite element (FEM) or finite difference (FDM) methods. For Helmholtz equations with high wavenumbers the discretized problem becomes extremely large, as the number of mesh points per wavenumber should be sufficiently large to result in an acceptable solution, which prohibits the use of direct methods. Iterative methods, such as Krylov subspace solvers, offer an interesting alternative. However, Krylov methods are not competitive when solving Helmholtz equations without a good preconditioner~\cite{Ernst2012}.\\

In a previous paper~\cite{YokotaArxiv13}, we proposed using the FMM as a preconditioner for Krylov subspace methods by equipping it with boundary integral capability for satisfying conditions at finite boundaries. Our model problems included inhomogeneous Poisson and Stokes equations and showed that the FMM preconditioner performs similarly to algebraic multigrid in convergence rate, while excelling in scalings. Multigrid methods, while enormously effective when applied to coercive equations, have severe convergence problems when applied to the indefinite Helmholtz equation~\cite{Trottenberg2001}, however. The reason for this is while the characteristic components of the Helmholtz problem can be accurately approximated by the discrete equations on the fine grids, these components are invisible to any local relaxation since their errors can have very small residuals. On the other hand, the characteristic components can not be approximated on coarser grids since these grids do not resolve their oscillations~\cite{Brandt1997}.

In this paper, we extend our previous work and employ the FMM as a preconditioner for the Helmholtz problem. Although FMM for the Laplace and Helmholtz equations share the same computational structure, the expansions they use are based on very different principles. The outline of this paper is as follows. In \sect \ref{sec:model} we present the Helmholtz model problem. Sections \ref{sec:krylov} and \ref{sec:pre} respectively give overview of Krylov subspace methods and preconditioning. The framework of our preconditioner is discussed in \sect \ref{sec:framework} and includes: FMM, the essential kernel that makes our method efficient and scalable, and the boundary element method (BEM), which is the basis of the FMM preconditioner. Our numerical results in \sect \ref{sec:results} examine the convergence rates of the geometric multigrid, algebraic multigrid, fast multipole, and incomplete Cholesky preconditioners for small 2D Helmholtz problems. In \sect \ref{sec:Performance} we scale up the problems and perform strong scalability runs on 3D model problems. Our conclusions are given in \sect \ref{sec:conc}.

\section{Mathematical model}
\label{sec:model}

The problem under consideration in this manuscript is the solution of the Helmholtz equation which is the frequency domain model of wave propagation. While most applications are concerned with waves propagation in exterior domains, it is common to utilize Helmholtz equations posed in interior domains with impedance boundary conditions to describe acoustic and elastic problems in finite domains. The Helmholtz equation takes the form
\begin{subequations}
\begin{alignat}{3}
\nabla^2 u + \kappa^2 u & = f && \text{ in } & \domain,\\
\partial_n u - i\kappa u & = g && \text{ on } & \bound, 
\end{alignat}
\label{eq:helm}
\end{subequations}
\noindent where $\domain$ is a connected bounded domain in $\R^d$ ($d = 2,3$) with piecewise smooth boundary $\bound$, $\kappa$ represents a constant wave number, and $f$ and $g$ are prescribed complex functions.

Discretizing~\eqref{eq:helm} by finite element or finite difference methods leads to a large sparse linear system of the form
\begin{equation}
A x = b, 
\end{equation}
where $A \in \Cnn$ is a large sparse symmetric matrix and $\bb \in \Cn$ contains the forcing and  boundary data. For large values of $\kappa$, the matrix $A$ is complex-valued and indefinite, i.e., $A$ has eigenvalues with both positive and negative real parts. 

Iterative methods, such as Krylov subspace solvers, are widely used in many areas of scientific computing for solving such large sparse linear systems where direct methods, although robust and reliable, have expensive computational requirements.

\section{Krylov subspace methods}
\label{sec:krylov}

The main idea of Krylov subspace methods is to generate a basis of Krylov subspace
\begin{equation}
\mathcal{K}_j(A,r_0)=\text{span}\{r_0,Ar_0,A^2r_0,\dots,A^{j-1}r_0\},
\end{equation}
and then seek an approximate solution to the original problem from this subspace. Here, $r_0 = b - A x_0$, $x_0$ is the initial approximate solution, and $\mathcal{K}_j(A,r_0)$ is the $j^{th}$ Krylov subspace associated with $A$ and $r_0$. A wide variety of iterative methods fall within the Krylov subspace framework. This section focuses on methods for solving linear systems with indefinite coefficient matrices. 

The Minimal Residual method (MINRES)~\cite{paige1975} can be used to solve linear systems with symmetric indefinite coefficient matrices, as well as its generalization to the nonsymmetric case, GMRES~\cite{Saad1986}. Both algorithms have the minimization property but GMRES has the advantage that theoretically it guarantees convergence. The main problem in GMRES is that it uses long recurrences which implies that the amount of storage increases at each iteration. Therefore, applications of GMRES may be limited by available storage. To overcome this problem, restarted versions of the GMRES method are used, e.g., GMRES($m$)~\cite{Saad1986}. In the restarted GMRES, computation and storage costs are limited by specifying a fixed number of vectors to be generated. However, since restarting removes the previous convergence history, GMRES($m$) does not guarantee convergence. The Bi-Conjugate Gradient Stabilized (BiCGSTAB)~\cite{vanderVorst1992} and Conjugate Gradient Squared (CGS)~\cite{Sonneveld1989} methods are short recurrence alternatives to GMRES. Even though BiCGSTAB is generally more stable and robust than CGS~\cite{broyden2004}, neither method guarantees monotonically decreasing residuals. We refer the reader to the books by Greenbaum~\cite{greenbaum1997} and Saad~\cite{saad2003} for more details on Krylov methods.

\begin{figure}
\centering
\includegraphics[width=0.5\textwidth, height=5cm]{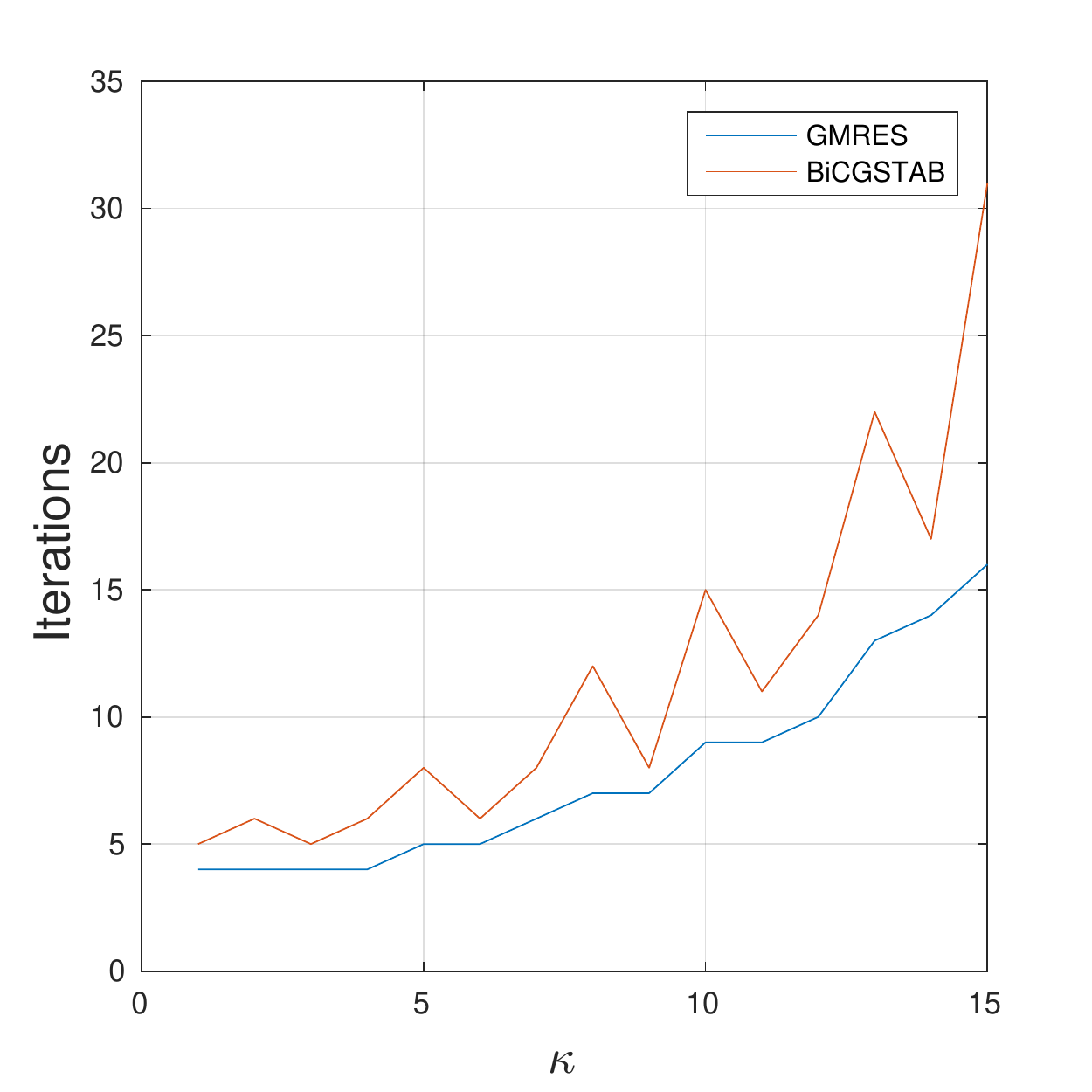}
\caption{Number of FMM-preconditioned GMRES and BiCGSTAB iterations required for the relative residual to reduce by six orders of magnitude for different values of $\kappa$ for the problem described in~\eqref{eq:helm1}, $h=2^{-5}$.}
\label{fig:gmres_bigcstab}
\end{figure}

The choice of an iterative method for indefinite problems is not straightforward. If the matrix-vector multiplication is expensive, e.g., when the coefficient matrix is dense, then GMRES is the method of choice since it requires the fewest matrix-vector multiplications to converge to the desired tolerance~\cite{greenbaum1997}. If the matrix-vector multiplication is not so expensive, then methods like BiCGSTAB and CGS are probably good choices. Figure~\ref{fig:gmres_bigcstab} compares, for different values of $\kappa$, the number of FMM-preconditioned GMRES and BiCGSTAB iterations for the relative residual to reduce by six orders of magnitude for the Helmholtz problem with homogeneous Dirichlet conditions on all boundaries and a nonhomogeneous source, as described in~\eqref{eq:helm1}. We do not consider CGS since sometimes its convergence curve shows wild oscillations that can lead to numerical instabilities~\cite{greenbaum1997}. Figure~\ref{fig:gmres_bigcstab} shows that the FMM-preconditioned GMRES outperformed BiCGSTAB for all values of $\kappa$. For this reason, we decided to use the GMRES iterative solver for all our test problems in Sections~\ref{sec:results} and~\ref{sec:Performance}.\\

\begin{figure}[t]
\centering
\subfigure[Laplace equation]{\includegraphics[width=0.35\textwidth]{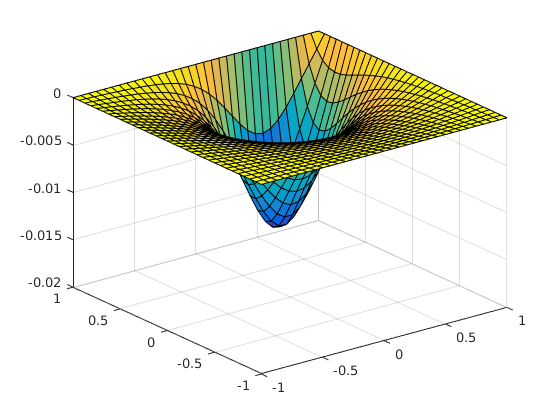}\label{f:k_0}}
\subfigure[Helmholtz equation]{\includegraphics[width=0.35\textwidth]{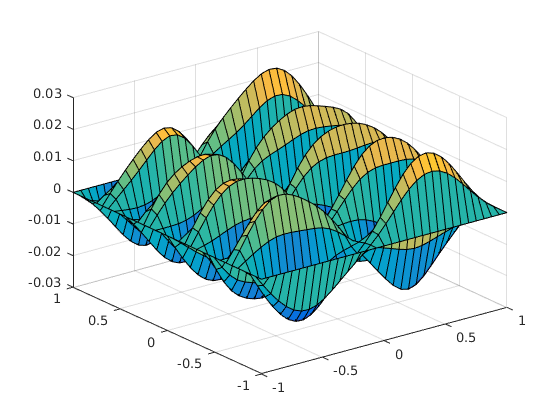}\label{f:k_15}}
\caption{Solution of Laplace and Helmholtz equations with the same boundary conditions.}
\label{f:k_0_15}
\end{figure}

\begin{figure}[t]
\centering
\subfigure[Laplace equation]{\includegraphics[width=0.4\textwidth,height=4cm]{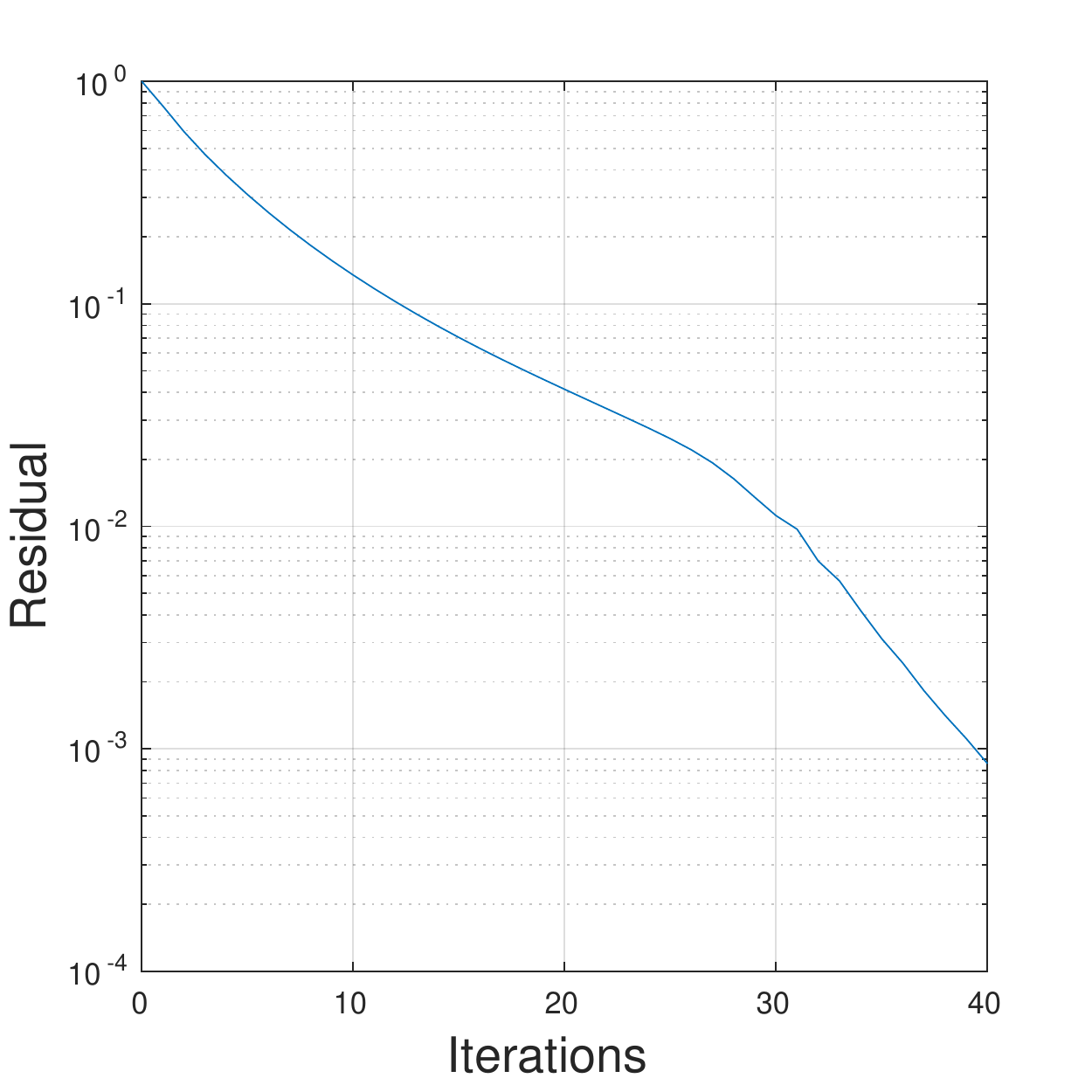}\label{f:k0}}
\subfigure[Helmholtz equation]{\includegraphics[width=0.4\textwidth,height=4cm]{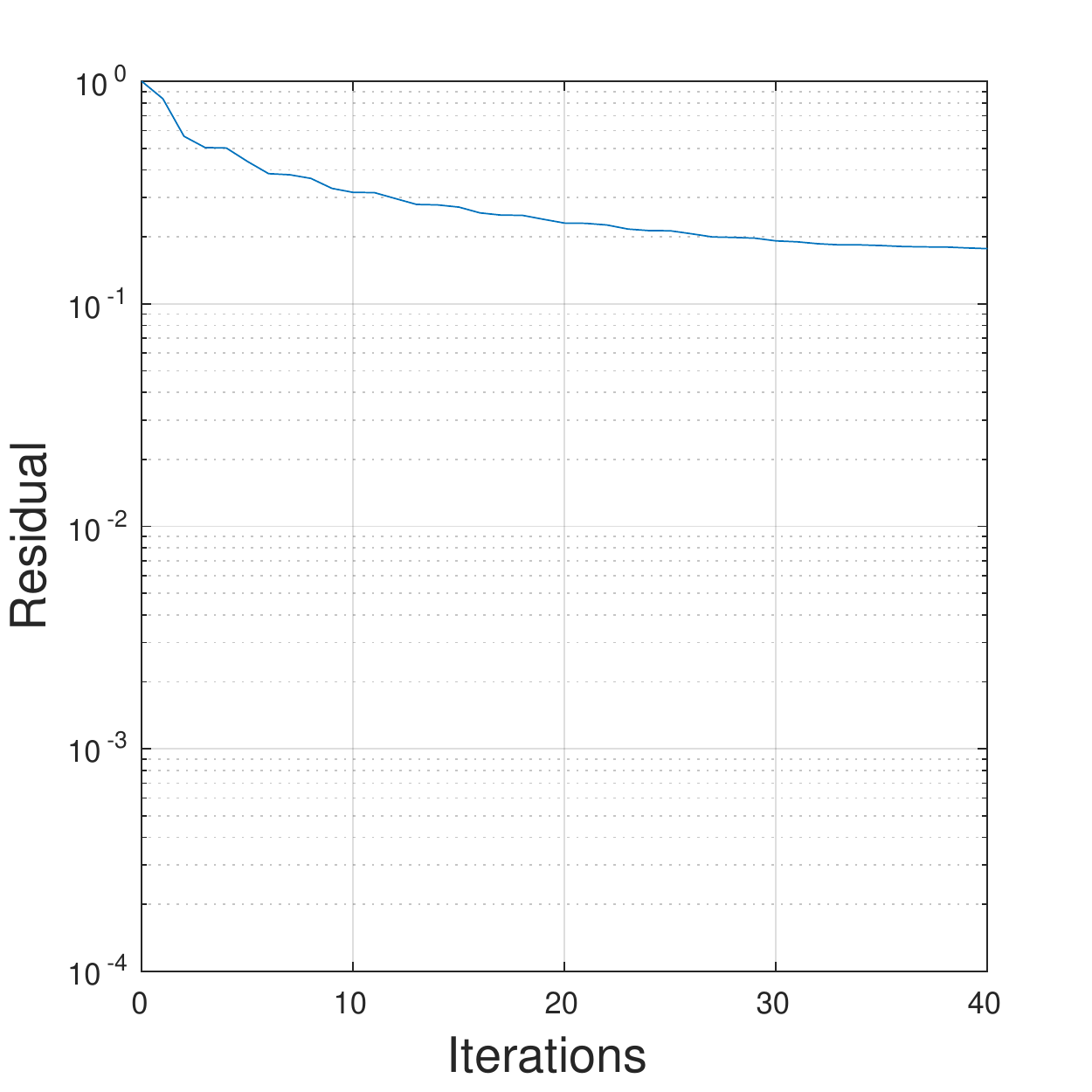}\label{f:k15}}
\label{fig:conv_k0_15}
\caption{Evolution of the residual of the unpreconditioned GMRES method for the Laplace equation, $\kappa = 0$, and the Helmholtz equation, $\kappa = 15$.
}
\end{figure}

Using the problem defined in~\eqref{eq:helm1}, Figure~\ref{f:k_0_15} shows the fundamental influence of the wave number $\kappa$ on the solution of the Helmholtz equation by comparing it against the Laplace equation ($\kappa = 0$). The solution of the Laplace equation is large only near the point source, Figure~\ref{f:k_0}, while for the Helmholtz equation, with $\kappa = 15$, the solution takes on large values periodically throughout the domain, Figure~\ref{f:k_15}. Figure~\ref{fig:conv_k0_15} shows how this influences the convergence of the unpreconditioned GMRES method. While the residual decreases rapidly for the Laplace equation, Figure~\ref{f:k0}, convergence stagnates for the Helmholtz equation, Figure~\ref{f:k15}. It is therefore important to have a preconditioner as GMRES method alone is not an effective iterative solver for Helmholtz equations~\cite{Ernst2012}.

\section{Preconditioning}
\label{sec:pre}

To improve the convergence of Krylov subspace methods, a preconditioner should be incorporated. The general rule for preconditioners is that the preconditioned system should be easy to solve, i.e., converges rapidly, and cheap to apply~\cite{Benzi2002}. It is important to strike a balance between these two requirements as they are competing with each other. 

One can apply a preconditioner on the left of the linear system, the right, or a combination of both. By left preconditioning, we solve a linear system premultiplied by a preconditioning matrix $M^{-1}$, i.e.,
\begin{equation}
M^{-1} A x = M^{-1} b.
\end{equation}
On the other hand, right preconditioning  is based on solving
\begin{equation}
A M^{-1} \hat{x} = b,
\end{equation}
\noindent where $\hat{x} = M x$. Both preconditionings show typically a similar convergence behavior and the type of preconditioning to use depends mainly on the choice of the iterative method and the problem characteristics. For example, right preconditioning is often used with the GMRES method~\cite{Benzi2002}. The essential difference between left and right preconditioned GMRES is that left-preconditioned GMRES computes residuals based on the preconditioned system while the residuals for the right-preconditioned GMRES are identical to the true residuals. This difference may affect the stopping criterion~\cite{saad2003}.\\

When Krylov subspace methods are used, it is not necessary to form the preconditioning matrix $M^{-1}$ explicitly. Instead, the preconditioning matrix can be a linear operation that defines the inverse of a matrix implicitly. This allows us to use matrix-free preconditioners such as the fast multipole method.

\section{Framework of the Fast Multipole preconditioner}
\label{sec:framework}

Fast Multipole methods have high arithmetic intensity, high degree of parallelism, controllable accuracy, and potentially less synchronous communication pattern compared to factorization based and multilevel methods. These features make the FMM a promising preconditioner for large systems on future architectures. One apparent disadvantage of the FMM approach to preconditioning is that the fast multipole method does not naturally incorporate boundary conditions, as these do not appear in its original application of solving $N$-body problems. Some approaches to overcome this obstacle include utilizing the boundary element method or the method of images. In the current work, we couple the FMM to the BEM since BEM gives the flexibility to solve over arbitrary geometries. This section briefly introduces some key ingredients of the FMM and BEM algorithms, more details on the FMM can be found in~\cite{greengard1987,Beatson1997}, and on the BEM in~\cite{sauter2011}.

\subsection{Overview of the Fast Multipole Method}
\label{subsec:fmm}

In the domain of scientific computing, $N$-body problems are used to simulate physical bodies or elementary particles interaction under physical forces that affect them from gravity or electromagnetic field~\cite{greengard1987}. $N$-body problem can be represented mathematically by the sum
\begin{equation}
f(y_j) = \sum_{i=1}^N w_i K(y_j,x_i),
\end{equation}
\noindent where $f(y_j)$ represents a field value evaluated at a point $y_j$ which is generated by the influence of sources located at the set of centers $\{x_i\}$, $\{x_i\}$ is the set of source points with weights given by $w_i$, $\{y_j\}$ is the set of evaluation points, and $K(y,x)$ is the kernel that governs the interactions between evaluation and source particles. 

The direct approach to simulate the $N$-body problem evaluates all pair-wise interactions among the particles and results in a computational complexity of $\mathcal{O}(N^2)$. This complexity is prohibitively expensive even for modestly large data sets. For simulations with large data sets, many faster algorithms have been invented, e.g., tree code~\cite{Barnes1986} and fast multipole methods~\cite{greengard1987}. The basic idea behind these fast algorithms is to cluster particles at successive levels of spatial refinement. The tree code clusters the far particles and achieves $\mathcal{O}(N \log N)$ complexity. The further apart the particles, the larger the interaction groups into which they are clustered. On the other hand, FMM divides the computational domain into near-domain and far-domain and computes interactions between clusters by means of local and multipole expansions, providing $\mathcal{O}(N)$ complexity.

\subsubsection{Hierarchical domain decomposition}
%
\begin{figure}
\centering
\subfigure[2D domain]{\includegraphics[width=0.2\textwidth]{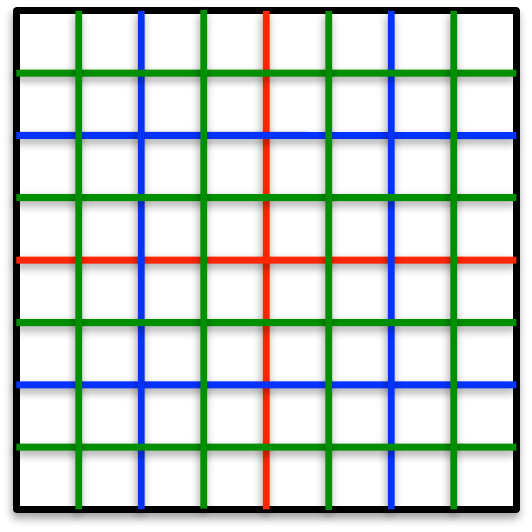}\label{f:rhs_1}\label{f:domain}}
\subfigure[Quad-tree]{\includegraphics[width=0.75\textwidth]{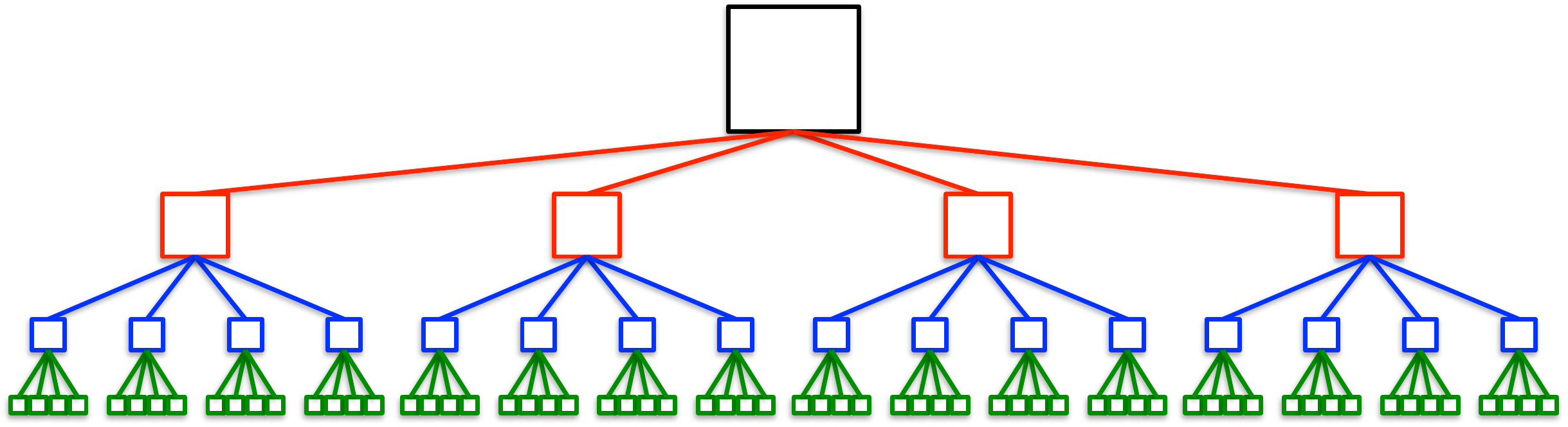}\label{f:tree}}
\caption{Decomposition of a 2D computational domain into a quad-tree.}
\label{fig:decomposition}
\end{figure}

The first step of the FMM algorithm is the decomposition of the computational domain. This spatial decomposition is accomplished by a hierarchical subdivision of the space associated to a tree structure. The 3D spatial domain of FMM is represented by oct-trees, where the space is recursively subdivided into eight cells until the finest level of refinement or ``leaf level". Figure~\ref{fig:decomposition} illustrates such hierarchical space decomposition for a 2D domain, Figure~\ref{f:domain}, associated to a quad-tree structure, Figure~\ref{f:tree}.

\subsubsection{The FMM Calculation Flow}

\begin{figure*}
\centering
\includegraphics[width=0.8\textwidth]{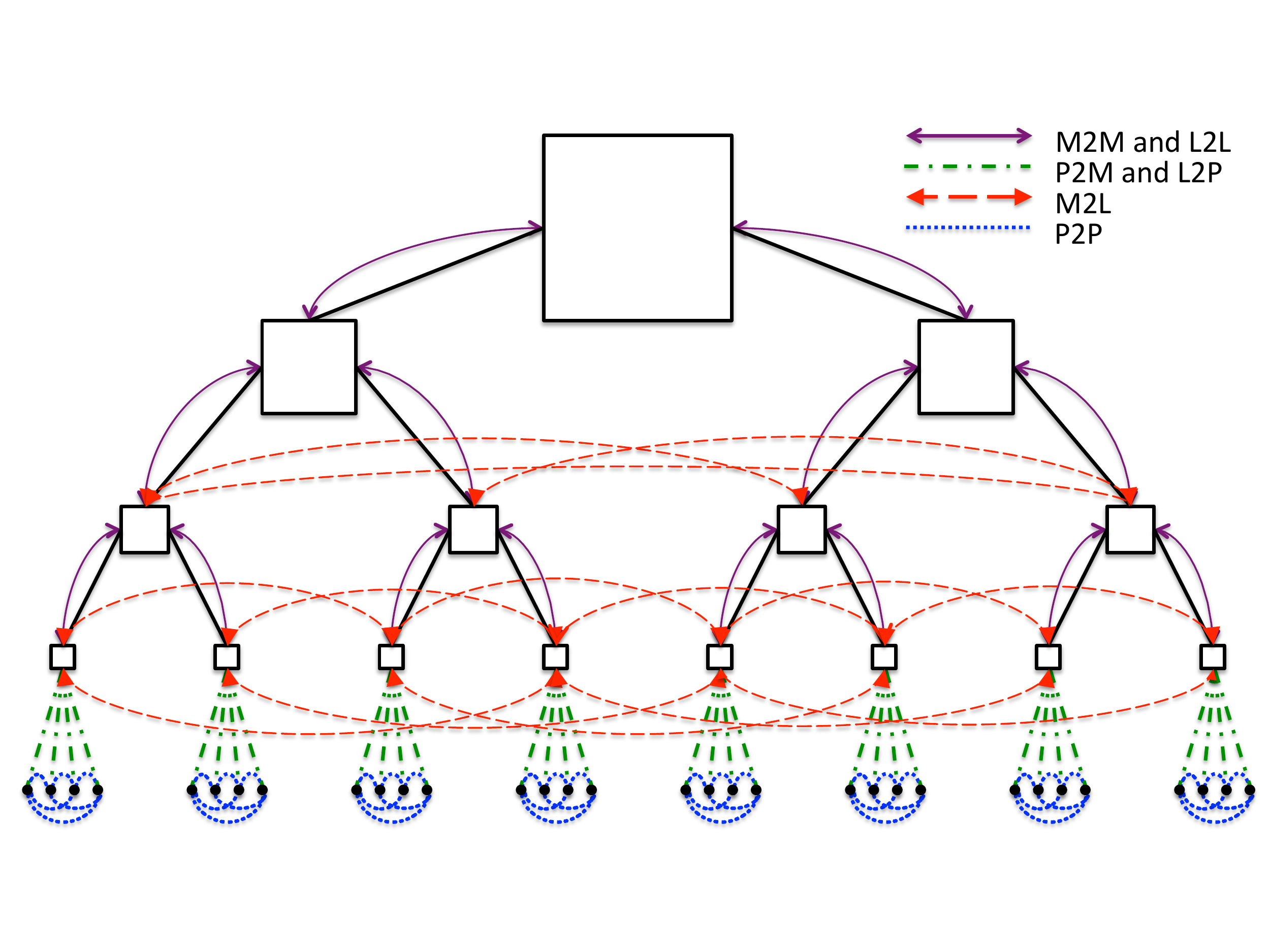}
\caption{Illustration of the FMM kernels: P2M (Particle-to-Multipole), M2M (Multipole-to-Multipole), M2L (Multipole-to-Local), L2L (Local-to-Local), L2P (Local-to-Particle), and P2P (Particle-to-Particle).}
\label{fig:fmm_flow}
\end{figure*}

The flow of the FMM calculation starts by transforming the mass/charge of the source particles into multipole expansions by means of a Particle-to-Multipole kernel (P2M). Then, the multipole expansions are translated to the center of larger cells using a Multipole-to-Multipole kernel (M2M). FMM calculates the influence of the multipoles on the target particles in three steps. First, it translates the multipole expansions to local expansions using a Multipole-to-Local kernel (M2L). Next, these local expansions are translated to smaller cells using a Local-to-Local kernel (L2L). Finally, the effect of local expansions in the far field is translated onto target particles using a Local-to-Particle kernel (L2P). All pairs interaction is used to calculate the effect of near field on the target particles by means of a Particle-to-Particle kernel (P2P). Figure~\ref{fig:fmm_flow} illustrates the FMM main kernels: Particle-to-Multipole (P2M), Multipole-to-Multipole (M2M), Multipole-to-Local (M2L), Local-to-Local (L2L), Local-to-Particle (L2P), and Particle-to-Particle (P2P).

\subsection{Conventional BEM formulation}
\label{subsec:bem}

In the conventional Galerkin boundary element method~\cite{sauter2011}, the boundary integral equations are obtained using Green's function, or Cauchy's integral theorem, and are then solved using discretization techniques similar to those employed in the finite element method (FEM). A brief description of the BEM formulation is given in this section. \\

Applying Green's theorem to~\eqref{eq:helm} gives the following Helmholtz integral equation
\begin{equation}
\eta(P) u(P)=\int_\Gamma\frac{\partial u}{\partial n}Gd\Gamma-\int_\Gamma u\frac{\partial G}{\partial n}d\Gamma+\int_\Omega fGd\Omega,
\label{eq:bem}
\end{equation}
where $P$ denotes a point in $\domain$, $n$ is the outward normal unit vector on $\bound$, and $\eta(P)$ is the Cauchy principal value of the boundary integral and is defined as
\[ \eta(P) =
  \begin{cases}
    \dfrac{1}{2}       & \quad \text{if } P \in \bound,\\
    1       & \quad \text{if } P \in \domain\backslash \bound,\\
    0       & \quad \text{if } P \notin \domain\cup \bound.
  \end{cases}
\]
The free space Green's function $G$ and its normal derivative in 2D are given by
\begin{subequations}
\begin{alignat}{2}
G(P,Q) &= \dfrac{i}{4} H_0^1(\kappa r), \\ 
\frac{\partial G}{\partial n}(P,Q) &= -\dfrac{i\kappa}{4} H_1^1(\kappa r)r_n,
\end{alignat}
\end{subequations}
\noindent where $r$ is the Euclidean distance between points $P$ and $Q$, and $H_0^1$ and $H_1^1$ are the Hankel function of the first kind, $0^{th}$ and $1^{st}$ order, respectively. For the 3D problem, the Green's function and its normal derivative are defined as
\begin{subequations}
\begin{alignat}{2}
G(P,Q) &= \dfrac{1}{4 \pi r} e^{i \kappa r}, \\ 
\frac{\partial G}{\partial n}(P,Q) &= \dfrac{i\kappa r - 1}{4 \pi r^2} e^{i \kappa r}.
\end{alignat}
\end{subequations}

To solve a given boundary value problem, we start by solving for values on the boundaries. In case of Dirichlet boundary conditions for example, the resulting equation can be written as
\begin{equation}
\int_\Gamma\frac{\partial u}{\partial n}Gd\Gamma=\int_\Gamma u\left(\frac{1}{2}\delta+\frac{\partial G}{\partial n}\right)d\Gamma-\int_\Omega fGd\Omega\quad \text{on}\ \partial\Omega, \label{eq:ubound}
\end{equation}
Here, all values on the right-hand side are known, and $\partial u/\partial n$ on the boundaries is obtained by solving~\eqref{eq:ubound}. Afterwards, the solution to the original Helmholtz equation $u$ can be obtained by solving the following equation
\begin{equation}
u=\int_\Gamma\frac{\partial u}{\partial n}Gd\Gamma-\int_\Gamma u\frac{\partial G}{\partial n}d\Gamma+\int_\Omega fGd\Omega\quad \text{in}\ \Omega. \label{eq:uinter}
\end{equation}
At this point, all values on the right-hand side are known so one can perform three integrations to obtain $u$ at the internal nodes, which is the solution to~\eqref{eq:helm}. 

\subsection{Discretization}

Discretization of all boundary/domain integrals in~\eqref{eq:ubound} and~\eqref{eq:uinter} is done in three steps:
\begin{enumerate}
\item Break the global integral into a sum of piecewise local integrals over each element. For example, the integration of the second term on the right-hand side of~\eqref{eq:bem} can be expressed as
\begin{equation}
\int_\Gamma u\frac{\partial G}{\partial n}d\Gamma=\sum_j\int_{\Gamma_j}u\frac{\partial G}{\partial n}d\Gamma_j.
\label{eq:piecewise}
\end{equation}
Where the piecewise integration is still performed analytically.

\item Break the local integral into the sum of contributions from the basis functions of each node that belongs to the element. Integration over a piecewise element $\Gamma_j$ can be obtained from
\begin{equation}
\int_{\Gamma_j}u\frac{\partial G}{\partial n}d\Gamma_j=|J_j|\sum_ku_{jk}\int_{-1}^1\phi_k(\xi)\frac{\partial G}{\partial n}d\xi,
\label{eq:basis}
\end{equation}
where the index $k$ sums over all nodes in the element, $J_j$ is the Jacobian of the $j^{th}$ element, and $\phi_k$ is the basis function of node $k$.

\item Integrate over each basis function using Gaussian quadratures. Equation~\eqref{eq:basis} is not completely discretized since it still requires analytical integration over the parametrized space $\xi$. In special cases where the basis function $\phi$ is of low order, and the Green's function $G$ is simple, the integration can be performed analytically. Though, a more general solution to this problem is provided through numerical integration using Gaussian quadratures
\begin{equation}
\int_{-1}^1\phi_k(\xi)\frac{\partial G}{\partial n}d\xi=\sum_l\phi_k(\xi_l)\frac{\partial G_{jkl}}{\partial n}w_l,
\label{eq:quadrature}
\end{equation}
where $\xi_l$ and $w_l$ are the parametrized coordinates and weights of the quadratures, respectively. The Green's function is a function of the location of the quadrature points, which depends on $j$, $k$, and $l$, and is therefore noted as $G_{jkl}$. Combining~\eqref{eq:piecewise},~\eqref{eq:basis}, and~\eqref{eq:quadrature} gives
\begin{equation}
\int_\Gamma u\frac{\partial G}{\partial n}d\Gamma=\sum_j|J_j|\sum_ku_{jk}\sum_l\phi_k(\xi_l)\frac{\partial G_{jkl}}{\partial n}w_l,
\label{eq:element}
\end{equation}
where the indices $j$, $k$, and $l$ correspond to the elements, nodes, and quadrature points, respectively.\\
\end{enumerate}

In the present work, we use constant elements so there are no nodal points at the corners of the square domain for the tests in Sections \ref{sec:results} and \ref{sec:Performance}. By applying the aforementioned discretization technique to all the integrals in~\eqref{eq:ubound}, we obtain the following matrix-vector representation
\begin{multline*}
N_\Gamma
\left\{
\phantom{
\begin{bmatrix}
\ddots\\
&G_{ij}\\
&&\ddots
\end{bmatrix}
}
\right.
\hspace{-24mm}
\overbrace{
\begin{bmatrix}
\ddots\\
&G_{ij}\\
&&\ddots
\end{bmatrix}
}^{N_\Gamma}
\underbrace{
\begin{bmatrix}
\vdots\\
\frac{\partial u_j}{\partial n}\\
\vdots
\end{bmatrix}
}_\text{unknown}
=
\overbrace{
\begin{bmatrix}
\hspace{-18mm}\ddots\\
\frac{1}{2}\delta_{ij}+\frac{\partial G_{ij}}{\partial n}\\
\hspace{18mm}\ddots
\end{bmatrix}
}^{N_\Gamma}
\begin{bmatrix}
\vdots\\
u_j\\
\vdots
\end{bmatrix}
-
\overbrace{
\begin{bmatrix}
\ddots\\
&G_{ij}\\
&&\ddots
\end{bmatrix}
}^{N_\Omega}
\begin{bmatrix}
\vdots\\
f_j\\
\vdots
\end{bmatrix},
\end{multline*}
where $N_\bound$ and $N_\Omega$ are the number of boundary and internal nodes, respectively. The diagonal term $G_{ii}$ is singular and can be determined analytically using the following formula in 2D~\cite{Peterson1998}
\begin{equation}
\int_{\Gamma_m}  H_0^1(\kappa r_m) d\Gamma_m = w_m+i\dfrac{2}{\pi}w_m\big[\ln(\dfrac{\gamma \kappa w_m}{4})-1)\big],
\end{equation}
where $r_m=0$ for the diagonal, $w_m$ is the width of the local element $\Gamma_m$, and $\gamma=1.781072418$ is the exponential of Euler's constant.

Similarly, we apply the same discretization technique to~\eqref{eq:uinter} and obtain the solution of the Helmholtz equation~\eqref{eq:helm}
\begin{multline*}
\small
N_\Omega
\left\{
\begin{bmatrix}
\vdots\\
u_i\\
\vdots
\end{bmatrix}
\right.
=
\overbrace{
\begin{bmatrix}
\ddots\\
&G_{ij}\\
&&\ddots
\end{bmatrix}
}^{N_\Gamma}
\begin{bmatrix}
\vdots\\
\frac{\partial u_j}{\partial n}\\
\vdots
\end{bmatrix}
-
\overbrace{
\begin{bmatrix}
\hspace{-18mm}\ddots\\
\frac{\partial G_{ij}}{\partial n}\\
\hspace{18mm}\ddots
\end{bmatrix}
}^{N_\Gamma}
\begin{bmatrix}
\vdots\\
u_j\\
\vdots
\end{bmatrix}
+
\overbrace{
\begin{bmatrix}
\ddots\\
&G_{ij}\\
&&\ddots
\end{bmatrix}
}^{N_\Omega}
\begin{bmatrix}
\vdots\\
f_j\\
\vdots
\end{bmatrix}.
\end{multline*}\\
The third term on the right-hand side is the dominant part of the computational load as it involves an $N_\Omega\times N_\Omega$ matrix. This matrix-vector multiplication can be approximated in $\mathcal{O}(N \log N)$ time using the FMM described in~\sect~\ref{subsec:fmm}. We also use the FMM to accelerate all the matrix-vector multiplications in the discretized forms of~\eqref{eq:ubound} and~\eqref{eq:uinter}.

\section{Numerical Results}
\label{sec:results}

In this section, several 2D numerical experiments are performed to assess the convergence of the FMM preconditioner for an increasing mesh size ($h^{-1}$) and/or increasing wavenumber ($\kappa$). The domain $\domain$ is discretized by finite element method in MATLAB using IFISS~\cite{elman2007,ifiss} where we construct the coefficient matrix by adding the stiffness matrix to $\kappa^2$ times the mass matrix. The discretization results in a large sparse symmetric linear system which we solve using the GMRES iterative solver with $maxit=20$. The stopping condition is based on the relative residual norm satisfying the tolerance $ < 10^{-6}$. If convergence has not been achieved after $maxit$ iterations, the computation is terminated; this is denoted by `---' in the results. We compare the fast multipole preconditioner against the incomplete Cholesky (IC) factorization~\cite{meijerink1977} implemented in MATLAB and the algebraic multigrid (AMG) and geometric multigrid (GMG) methods from IFISS. For all problems and preconditioners the initial iterate is the zero vector. Our MATLAB implementation of the FMM used in this section is a direct $N$-body summation that is subsequently degraded to simulate a truncated FMM.\\

Our first Helmholtz problem~\cite{Erlangga2003} is posed on a unit square $[0,1]^2$ of homogeneous medium with homogeneous Dirichlet boundary conditions as follows:
\begin{subequations}
\begin{alignat}{3}
\nabla^2 u + \kappa^2 u & = (\kappa^2 - 5 \pi ^2) \sin(\pi x)\sin(2 \pi  y)&& \text{ in } & \domain,\label{eq:helm2_1}\\
u & = 0 && \text{ on } & \bound.
\end{alignat}
\label{eq:helm_fixed}
\end{subequations}
The exact solution of~\eqref{eq:helm_fixed} is
\begin{equation}
u = \sin(\pi x)\sin(2 \pi  y).
\end{equation}

\begin{table}
\centering
\caption{Effect of different FEM discretizations on the number of preconditioned GMRES iterations required for the relative residual to reduce by six orders of magnitude for the problem in~\eqref{eq:helm_fixed}, $\kappa=15$.}
\begin{tabular}{ccccc}
\hline\noalign{\smallskip}
          & \multicolumn{2}{c}{$Q_1$} & \multicolumn{2}{c}{$Q_2$} \\ 
          \hline\noalign{\smallskip}
 $h$      & AMG & FMM & AMG & FMM  \\
 \noalign{\smallskip}\hline\noalign{\smallskip}
 $2^{-5}$ &  14 & 4 & 8  & 6 \\ 
 $2^{-6}$ &  9  & 4 & 9  & 6 \\ 
 $2^{-7}$ &  9  & 3 & 9  & 5 \\
\noalign{\smallskip}\hline
\end{tabular}
\label{t:p4}
\end{table}

First, we study the effect of various finite element discretizations on the convergence of the FMM preconditioner. In particular, we consider the two FEM discretizations available from IFISS: $Q_1$, piecewise linear elements, and $Q_2$, piecewise quadratic elements. Table~\ref{t:p4} compares the number of GMRES iterations required for convergence to six digits of accuracy for the AMG and FMM preconditioners when using the $Q_1$ and $Q_2$ discretizations. While changing the finite element discretization has small effect on the AMG preconditioner for this particular problem, the number of FMM preconditioned GMRES iterations varies for different discretizations. Both $Q_1$ and $Q_2$ show mesh-independent convergence but the $Q_1$ discretization results in a smaller number of iterations for all the given values of $h$. For this reason, we use the $Q_1$ discretization for all test problems in this section. 

\begin{figure}
\centering
\subfigure[$\kappa = 5$]{\includegraphics[width=0.49\textwidth,height=2cm]{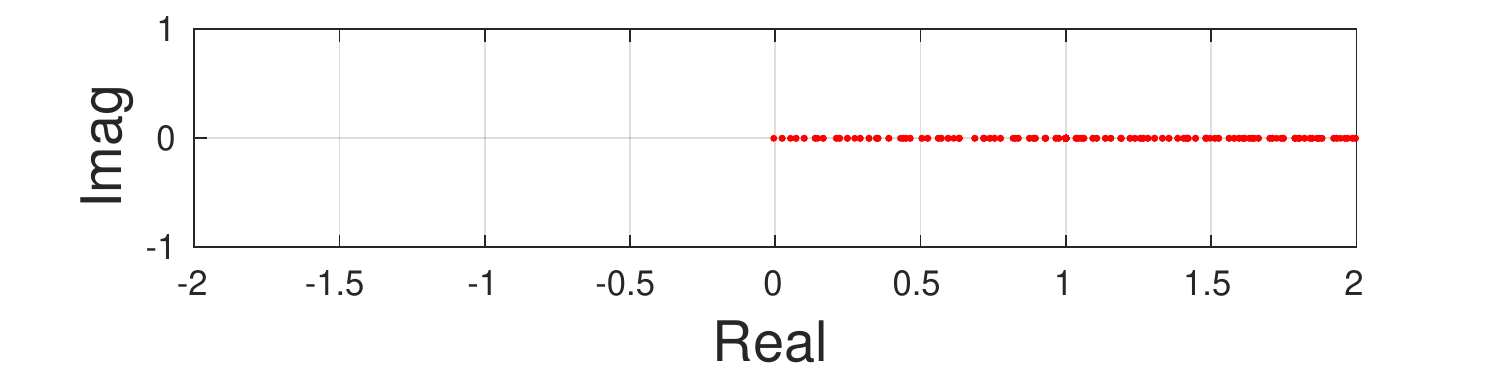}}
\subfigure[$\kappa = 10$]{\includegraphics[width=0.49\textwidth,height=2cm]{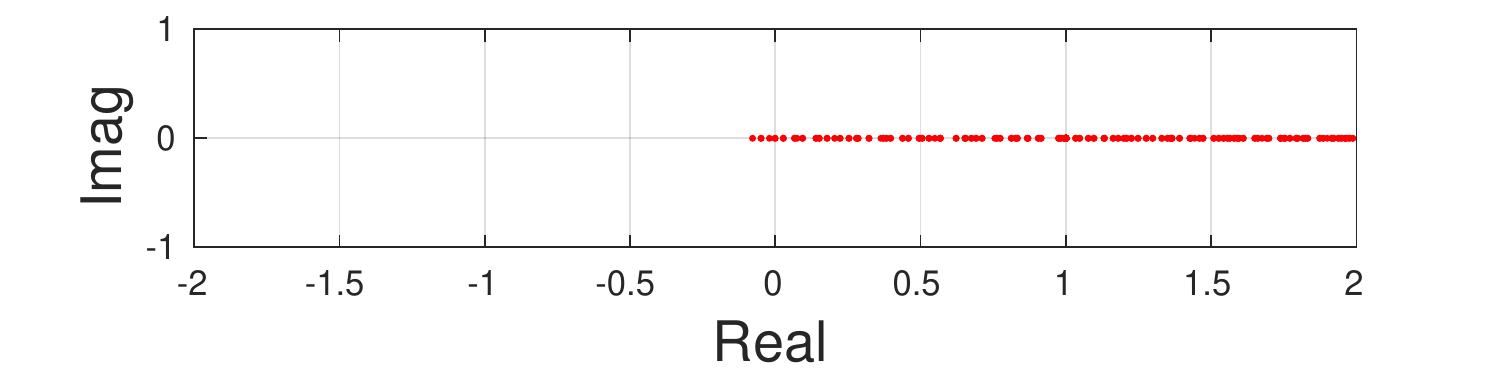}}
\subfigure[$\kappa = 20$]{\includegraphics[width=0.49\textwidth,height=2cm]{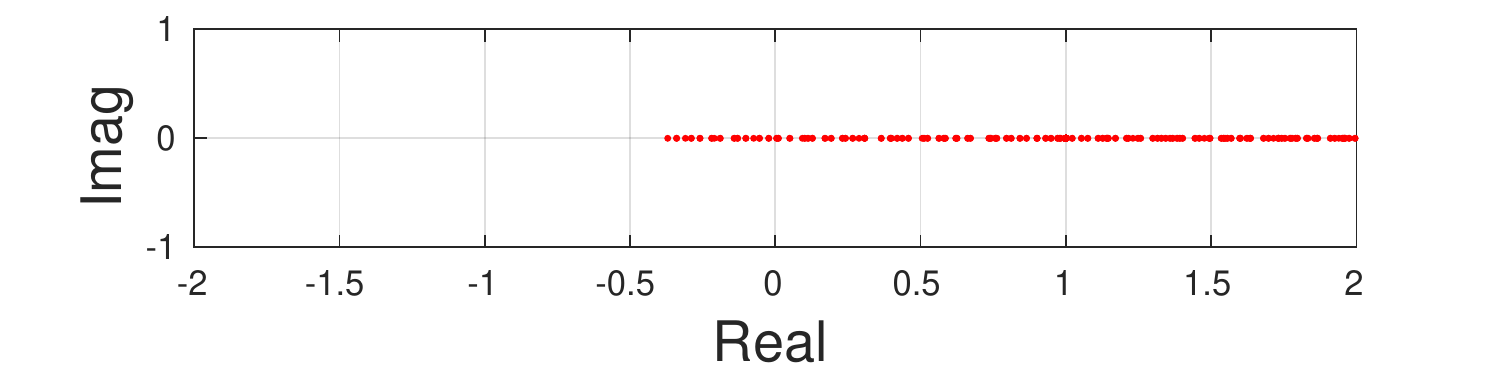}}
\subfigure[$\kappa = 40$]
{\includegraphics[width=0.49\textwidth,height=2cm]{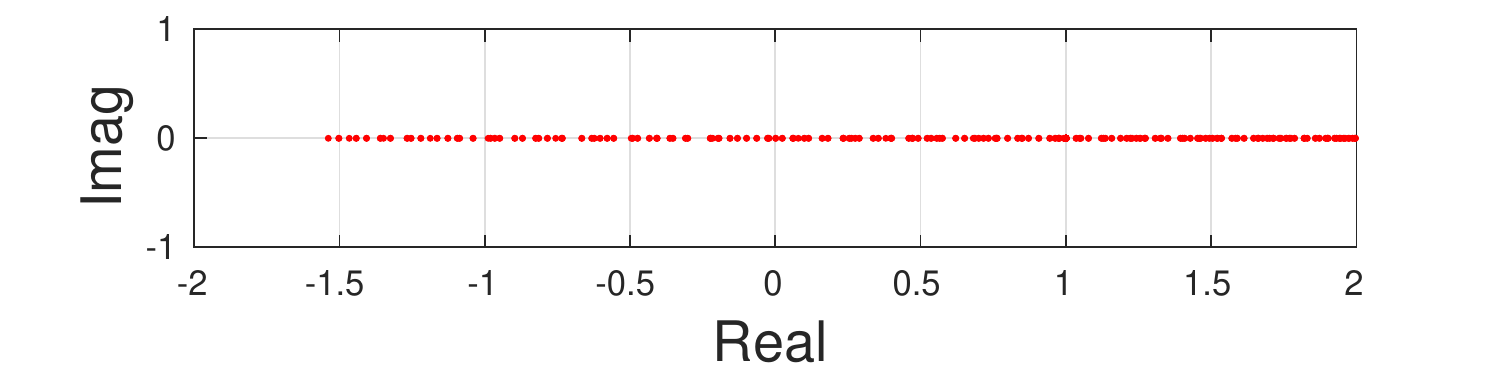}}
\caption{Eigenvalues of the the coefficient matrix $A$ for the problem descriped in~\eqref{eq:helm_fixed} with $\kappa$ = 5, 10, 20, and 40, respectively, $h=2^{-5}$.}
\label{fig:eigen_fixed}
\end{figure}

\begin{table}
\centering
\caption{Preconditioned GMRES iterations for the relative residual to reduce by six orders of magnitude for the problem in~\eqref{eq:helm_fixed}, $\kappa h=0.3125$.}
\label{t:helm_fixed}
\begin{tabular}{cccccc}
\hline\noalign{\smallskip}
 $h$      & $\kappa$  & GMG & AMG  & FMM &  IC   \\ 
\noalign{\smallskip}\hline\noalign{\smallskip}
 $2^{-4}$ & $5$       &  11 &   5  &  4  &  13 \\ 
 $2^{-5}$ & $10$      & --- &   6  &  4  &  --- \\ 
 $2^{-6}$ & $20$      & --- &  --- &  4  &  --- \\ 
 $2^{-7}$ & $40$      & --- &  --- &  5  &  --- \\
\noalign{\smallskip}\hline
\end{tabular}
\end{table}

The sub-figures in Figure~\ref{fig:eigen_fixed} show the eigenvalues of the coefficient matrix of the original linear system for $\kappa$ = 5, 10, 20, and 40, respectively. Notice that the coefficient matrix becomes more indefinite as the wavenumber increases with $\kappa=5$ resembles the definite problem. Different mesh refinements are used to solve~\eqref{eq:helm_fixed} with various wave numbers $\kappa$. Table~\ref{t:helm_fixed} lists the number of GMRES iterations required to reach the pre-specified tolerance for the GMG, AMG, FMM, and IC preconditioners. For small wave numbers, all preconditioners show a satisfactorily performance. GMG and IC become less effective for increasing values of $\kappa$ where the number of iterations required for convergence increases rapidly. For larger $\kappa$, preconditioning with FMM shows the best performance where other preconditioners fail to converge within $maxit$.\\

\begin{figure}
\centering
\subfigure[$f=e^{-10((y-1)^2 + (x-0.5)^2)}$]{\includegraphics[width=0.35\textwidth]{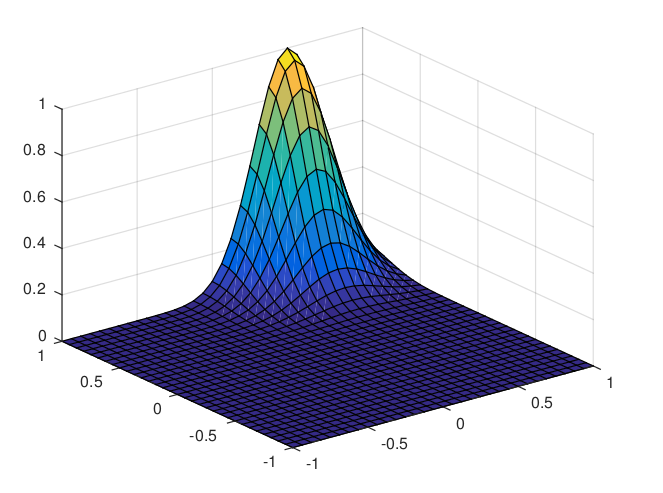}\label{f:rhs_1}}
\subfigure[Solution $u$ with $\kappa=5$]{\includegraphics[width=0.35\textwidth]{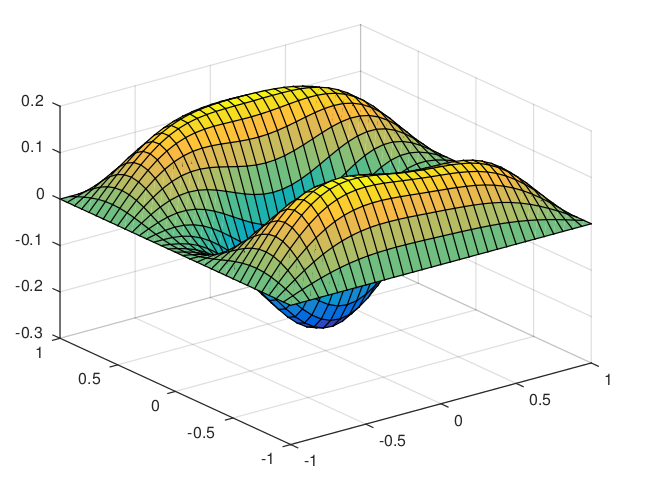}\label{f:direct_interaction}}
\caption{The right-hand term and solution of~\eqref{eq:helm1}.}
\label{fig:rhs_sol}
\end{figure}

The second Helmholtz example is posed on the square domain $[-1,1]^2$ and is characterized by homogeneous Dirichlet boundary conditions as follows:
\begin{subequations}
\begin{alignat}{3}
\nabla^2 u + \kappa^2 u & = e^{-10((y-1)^2 + (x-0.5)^2)} && \text{ in } & \domain,\label{eq:helm1_1}\\
u & = 0 && \text{ on } & \bound.
\end{alignat}
\label{eq:helm1}
\end{subequations} 
The right-hand side of~\eqref{eq:helm1_1} and the solution of~\eqref{eq:helm1} are shown in Figure~\ref{fig:rhs_sol}.
%

\begin{table}
\centering
\caption{Preconditioned GMRES iterations for the relative residual to reduce by six orders of magnitude for the problem in~\eqref{eq:helm1}, $\kappa=5$.}
\label{t:p1}
\begin{tabular}{ccccc}
\hline\noalign{\smallskip}
 $h$ & GMG & AMG & FMM & IC   \\ 
\noalign{\smallskip}\hline\noalign{\smallskip}
 $2^{-5}$ &  --- &  15  &  5  &  ---\\ 
 $2^{-6}$ &  --- &  15  &  4  &  ---\\ 
 $2^{-7}$ &  --- &  15  &  4  &  ---\\
\noalign{\smallskip}\hline
\end{tabular}
\end{table}

\begin{table}
\centering
\caption{Preconditioned GMRES iterations for the relative residual to reduce by six orders of magnitude for the problem in~\eqref{eq:helm1}, $h=2^{-6}$.}
\label{t:p2}
\begin{tabular}{ccccc}
\hline\noalign{\smallskip}
 $\kappa$ & GMG & AMG  & FMM & IC   \\ 
\noalign{\smallskip}\hline\noalign{\smallskip}
 $0.8$    &  5  &   5  &  4  &  --- \\ 
 $2$      &  8  &   6  &  4  &  --- \\ 
 $5$      & --- &  15  &  4  &  --- \\
\noalign{\smallskip}\hline
\end{tabular}
\end{table}

Table~\ref{t:p1} gives the number of GMRES iterations required for convergence on various mesh sizes with $\kappa=5$ for the GMG, AMG, FMM and IC preconditioners. Both FMM and AMG preconditioners appear to give mesh independent convergence, whereas GMG and IC factorization fail to converge within $maxit$. Table~\ref{t:p2} lists the number of preconditioned GMRES iterations for each preconditioner applied as a function of the wave number $\kappa$ with $h=2^{-6}$. For small wave numbers, the GMG, AMG, and FMM preconditioners show a very comparable performance. As $\kappa$ increases, both multigrid methods start to diverge while the FMM preconditioner maintains a wavenumber-independent convergence for the given range of $\kappa$.\\

The third Helmholtz example~\cite{Hui2007} is posed on $[0,1]^2$ and is characterized by inhomogeneous Dirichlet boundary conditions as follows:
\begin{subequations}
\begin{alignat}{3}
\nabla^2 u + \kappa^2 u & = 2 \sin (\mu x) \cos(\mu y)+4\mu x \cos(\mu x) \cos(\mu y) && \text{ in } & \domain,\label{eq:helm2_1}\\
u & = x^2 \sin(\mu x) \cos(\mu y) && \text{ on } & \bound,
\end{alignat}
\label{eq:helm2}
\end{subequations}
\noindent where $\kappa=\mu \sqrt{2}$.

\begin{figure}
\centering
\subfigure[$f$, $\mu=1$]{\includegraphics[width=0.24\textwidth]{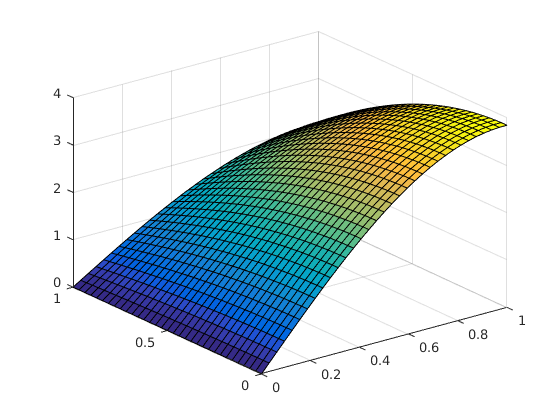}\label{f:direct_interaction}}
\subfigure[$u$, $\mu=1$]{\includegraphics[width=0.24\textwidth]{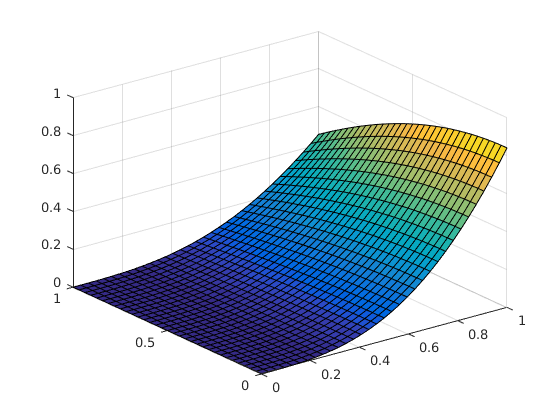}\label{f:direct_interaction}}
\subfigure[$f$, $\mu=4$]{\includegraphics[width=0.24\textwidth]{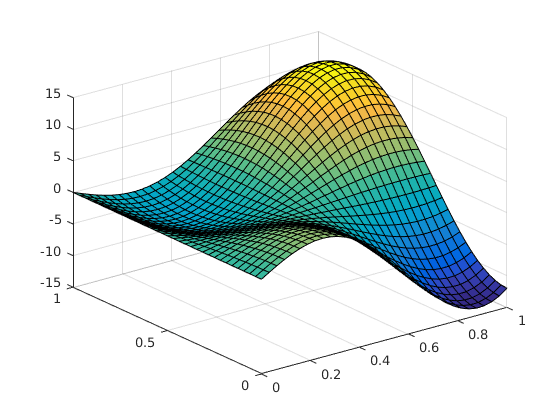}\label{f:fmm_interaction}}
\subfigure[$u$, $\mu=4$]
{\includegraphics[width=0.24\textwidth]{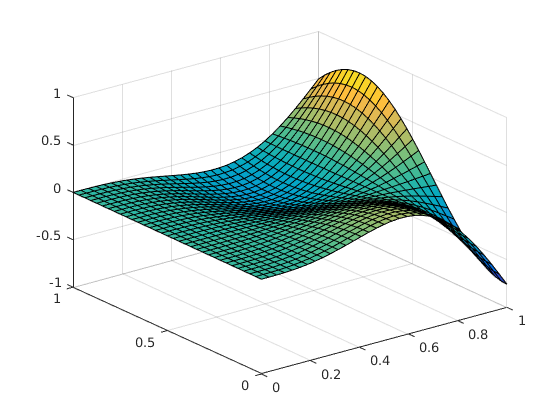}\label{f:direct_interaction}}\\
\subfigure[$f$, $\mu=6$]{\includegraphics[width=0.24\textwidth]{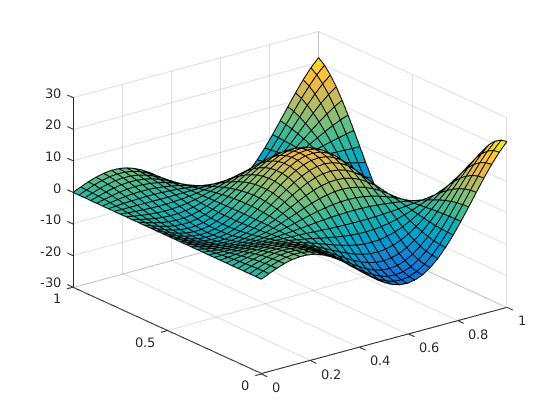}\label{f:fmm_flow}}
\subfigure[$u$, $\mu=6$]
{\includegraphics[width=0.24\textwidth]{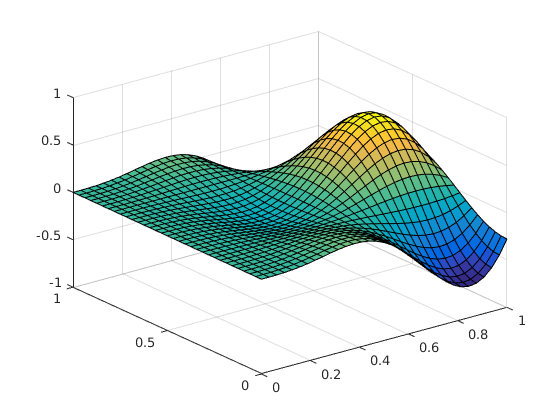}\label{f:direct_interaction}}
\subfigure[$f$, $\mu=8$]{\includegraphics[width=0.24\textwidth]{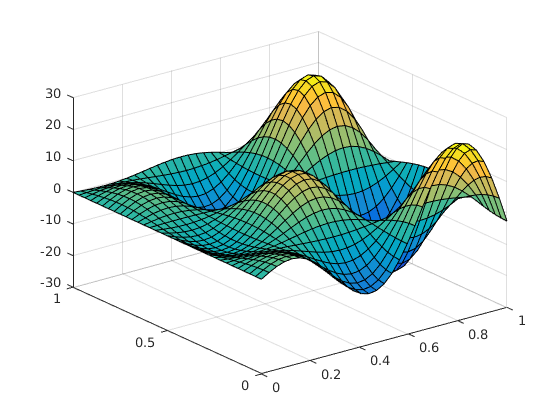}\label{f:fmm_flow}}
\subfigure[$u$, $\mu=8$]
{\includegraphics[width=0.24\textwidth]{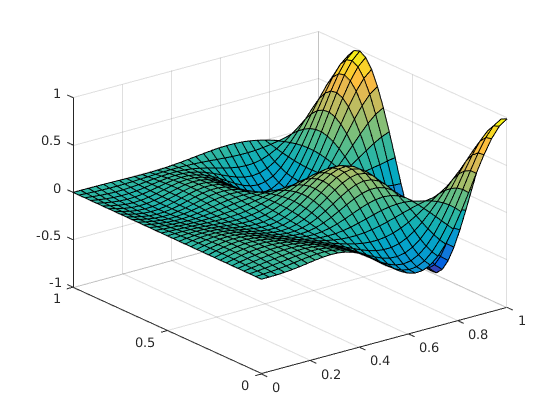}\label{f:direct_interaction}}\\
\caption{The right-hand term $f$ and solution $u$ of~\eqref{eq:helm2} for $\mu$ = 1, 4, 6, and 8.}
\label{fig:helm2}
\end{figure}

\begin{table}
  \centering
  \caption{Preconditioned GMRES iterations for the relative residual to reduce by six orders of magnitude for the problem in~\eqref{eq:helm2}.}
  \subtable[$\mu=1$]{
\begin{tabular}{ccccc}
\hline\noalign{\smallskip}
$h$ & GMG & AMG & FMM & IC\\ 
\noalign{\smallskip}\hline\noalign{\smallskip}
 $2^{-5}$ &  11  &  5  &  7  &  19  \\ 
 $2^{-6}$ &  14  &  5  &  6  &  ---  \\ 
 $2^{-7}$ &  17  &  5  &  6  &  ---  \\
\noalign{\smallskip}\hline
\end{tabular}}
\hspace{0.35 cm}
\subtable[$\mu=4$]{
\begin{tabular}{ccccc}
\hline\noalign{\smallskip}
$h$ & GMG & AMG & FMM & IC\\ 
\noalign{\smallskip}\hline\noalign{\smallskip}
 $2^{-5}$ &   16 &  7  &  7  &  ---  \\ 
 $2^{-6}$ &  --- &  7  &  6  &  ---  \\ 
 $2^{-7}$ &  --- &  7  &  6  &  ---  \\
\noalign{\smallskip}\hline
\end{tabular}}\\
  \subtable[$\mu=6$]{
\begin{tabular}{ccccc}
\hline\noalign{\smallskip}
$h$ & GMG & AMG & FMM & IC\\ 
\noalign{\smallskip}\hline\noalign{\smallskip}
 $2^{-5}$ &  --- &  8  &  6  &  ---  \\ 
 $2^{-6}$ &  --- &  8  &  6  &  ---  \\ 
 $2^{-7}$ &  --- &  8  &  6  &  ---  \\
\noalign{\smallskip}\hline
\end{tabular}}
\hspace{0.35 cm}
\subtable[$\mu=8$]{
\begin{tabular}{ccccc}
\hline\noalign{\smallskip}
$h$ & GMG & AMG & FMM & IC\\ 
\noalign{\smallskip}\hline\noalign{\smallskip}
 $2^{-5}$ &  --- &  15 &  7  &  ---  \\ 
 $2^{-6}$ &  --- &  15 &  7  &  ---  \\ 
 $2^{-7}$ &  --- &  15 &  6  &  ---  \\ 
\noalign{\smallskip}\hline
\end{tabular}}
\label{t:p3}
\end{table}

The right-hand term $f$ and exact solution $u$ in a unit square domain with various parameter $\mu$ are shown in Figure~\ref{fig:helm2}. Notice that the larger the value of $\mu$, the greater the variation of the function. The subtables in Table~\ref{t:p3} list the number of preconditioned GMRES iterations for the GMG, AMG, FMM, and IC preconditioners for $\mu$ = 1, 4, 6, and 8, respectively. For low frequencies, the GMG, AMG, and FMM preconditioners show a very satisfactorily comparable performance. The GMG preconditioner becomes less effective for increasing values of $\kappa$ where the number of iterations required for convergence exceeds $maxit$. For larger $\kappa$, the FMM preconditioner requires the smallest number of iterations to converge to the predefined tolerance. Examining all the subtables in Table~\ref{t:p3} shows that the FMM preconditioner achieves both mesh-independent and wavenumber-independent convergence for the given values of $h$ and $\kappa$.

\subsection{Impact of the FMM accuracy}

Fast multipole method has tunable accuracy, through control of the order of expansion, that enables it to trade-off accuracy for speed. In practice, this accuracy can reach machine precision if needed. When using FMM as a preconditioner for Krylov solvers, the FMM accuracy should be considered along with the accuracy of the BEM and the stopping criterion of the iterative solver. In our test problems, the order of expansion for the FMM is set to $p = 6$ which gives about six significant digits of accuracy. However, since we are using FMM as a preconditioner, the accuracy requirements are somewhat lower than that of general applications of FMM. This allows us to use low-accuracy FMM which, in practice, is significantly faster than the high-accuracy FMM even if it requires a few more iterations.

\begin{figure}
\centering
\includegraphics[width=0.7\textwidth,height=7cm]{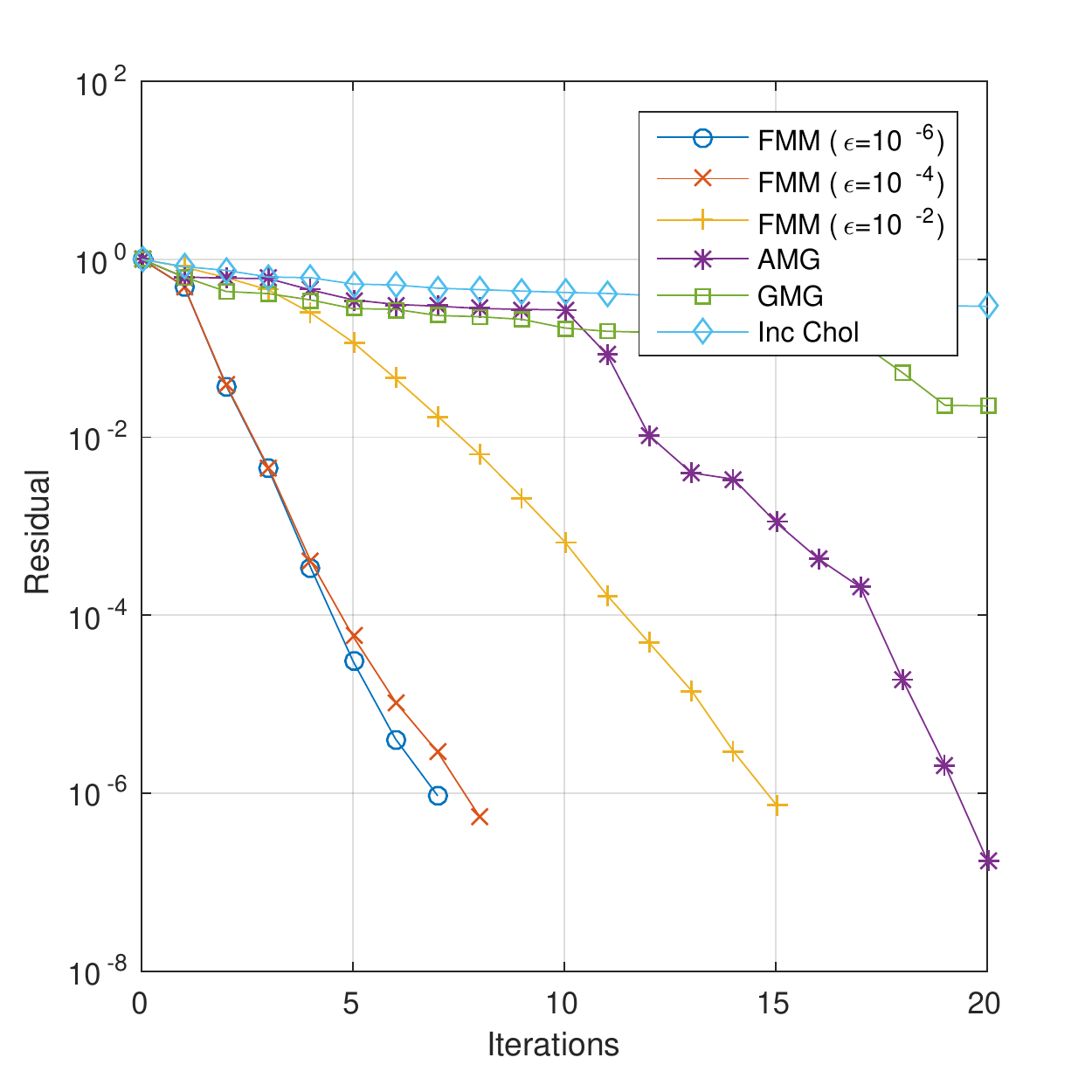}
\caption{Convergence rate of the FMM preconditioner with different precision, plotted along with AMG, GMG, and IC preconditioners. The $\epsilon$ represents the precision of the FMM where $\epsilon=10^{-6}$ corresponds to six significant digits of accuracy, $h=2^{-5}$.}
\label{fig:convergence}
\end{figure}

In Figure~\ref{fig:convergence} the relative residual at each GMRES iteration is plotted against the number of iterations for the fast multipole, algebraic multigrid, geometric multigrid, and incomplete Cholesky preconditioners using the problem defined in~\eqref{eq:helm1} with $\kappa = 7$. Three cases of the FMM are considered with six, four, and two significant digits of accuracy where $\epsilon=10^{-6}$ in the figure corresponds to the condition for the previous test problems. Decreasing the FMM accuracy to four digits slows down the convergence slightly. Decreasing the accuracy further to two digits slows down the convergence significantly, but is still better than the AMG, GMG and IC preconditioners. \\

\begin{figure}
\centering
\subfigure[$\lambda(A)$]{\includegraphics[width=0.49\textwidth,height=2cm]{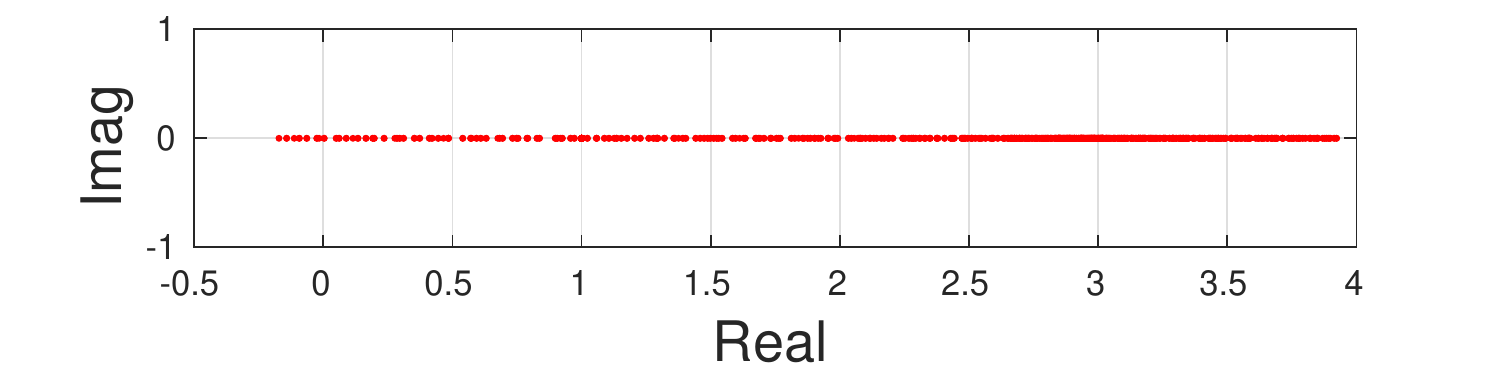}}
\subfigure[$\lambda(M^{-1}A)$, $\epsilon=10^{-2}$]{\includegraphics[width=0.49\textwidth,height=2cm]{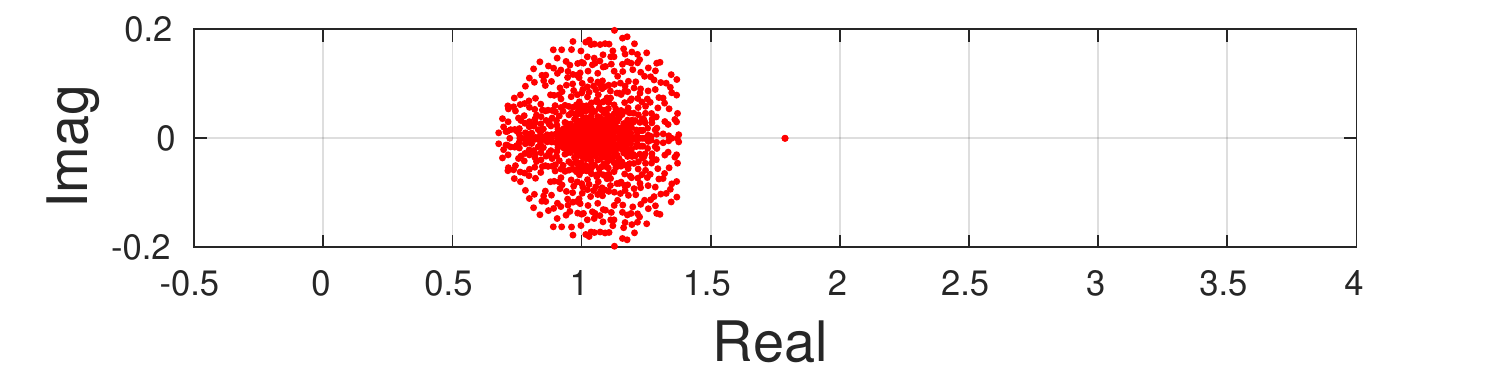}}
\subfigure[$\lambda(M^{-1}A)$, $\epsilon=10^{-4}$]{\includegraphics[width=0.49\textwidth,height=2cm]{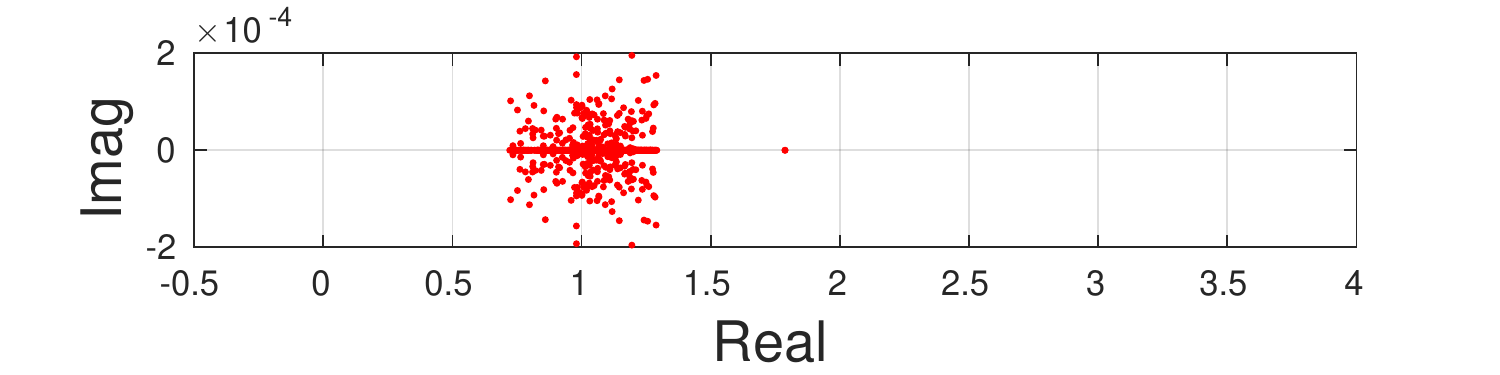}}
\subfigure[$\lambda(M^{-1}A)$, $\epsilon=10^{-6}$]
{\includegraphics[width=0.49\textwidth,height=2cm]{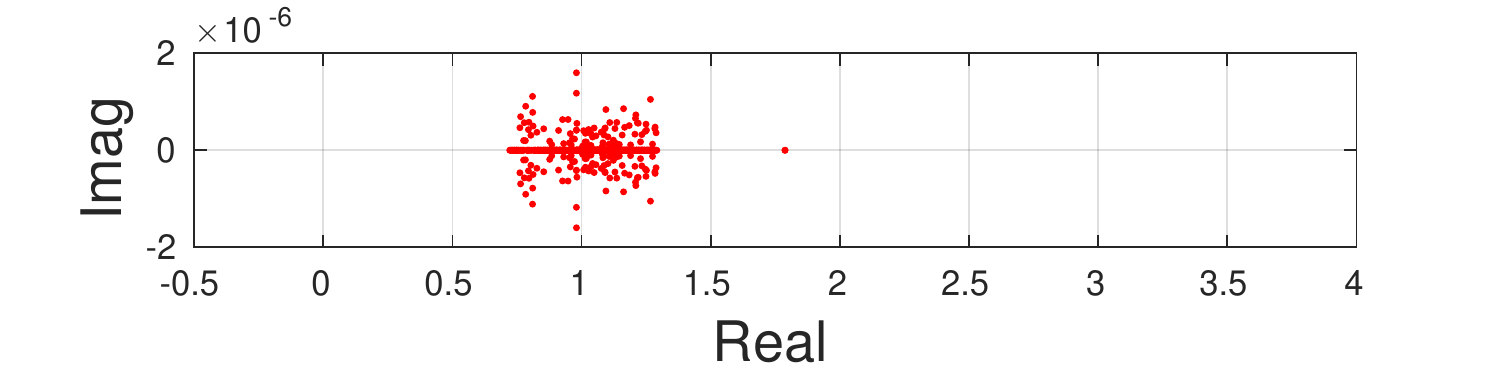}}
\caption{Eigenvalues of the coefficient matrix $A$ and the FMM-preconditioned matrix $M^{-1}A$ with different FMM precisions $\epsilon$ = $10^{-2}$, $10^{-4}$, and $10^{-6}$, $h=2^{-5}$.}
\label{fig:eigen}
\end{figure}

The best preconditioners for Krylov subspace methods move the smallest eigenvalues of the coefficient matrix away from the origin as small eigenvalues are known to often hamper the performance of Krylov solvers. The sub-figures in Figure~\ref{fig:eigen} show the eigenvalues clustering for the FMM-preconditioned coefficient matrix with different FMM precisions $\epsilon=$ $10^{-2}$, $10^{-4}$, and $10^{-6}$, respectively. Notice that as the FMM accuracy increases, the eigenvalues of the FMM-preconditioned matrix become better clustered and more bounded away from zero.

\section{Performance analysis}
\label{sec:Performance}

In this section, we evaluate the performance and scalability of the FMM-based preconditioner by implementing it in PETSc~\cite{petsc-user-ref,petsc-web-page}, via PetIGA~\cite{PetIGA} which is a software framework that sits on top of PETSc and facilitates a NURBS-based Galerkin finite element method, popularly known as isogeometric analysis (IGA).

All calculations were performed on Shaheen \RNum{2} which is a Cray XC40 with 6174 compute nodes, each with two 16-core Intel Haswell CPUs (Intel\textregistered Xeon\textregistered E5-2698 v3). The nodes of Shaheen \RNum{2} are connected by a dragonfly network using the Aries interconnect where the routers in each group are arranged as rows and columns of a rectangle with all-to-all links across each row and column but not diagonally. We use the GNU compiler and configured PETSc with  ``\texttt{COPTFLAGS=-O3 FOPTFLAGS=-O3 --with-clanguage=cxx\\ --download-fblas-lapack --download-hypre --download-metis \\--download-parmetis --download-superlu\_dist --with-debugging=0}". \\
All codes used in this work are made publicly available. A branch of PetIGA that includes the FMM preconditioner is hosted on Bitbucket~\footnote{https://bitbucket.org/rioyokota/petiga-fmm} and the open source FMM library ExaFMM is available on Github.~\footnote{https://github.com/exafmm/exafmm.git}\\

\begin{figure}[t]
\centering
\includegraphics[width=0.7\textwidth,height=7cm]{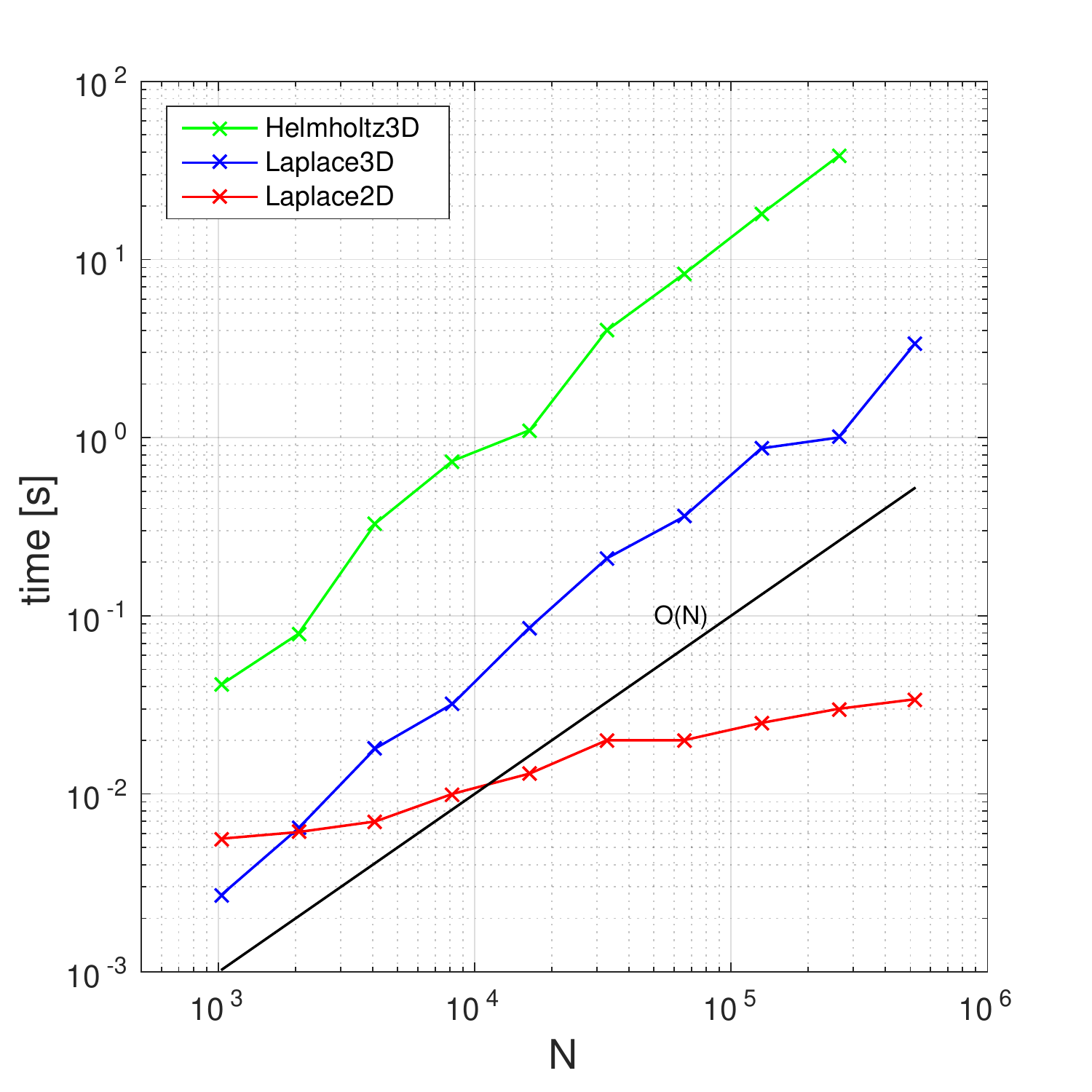}
\caption{Calculation time of the Poisson and Helmholtz FMM for the same  problem size.}
\label{fig:fmm_kernels}
\end{figure}

In a previous work~\cite{YokotaArxiv13}, we compared the strong scalability of the FMM preconditioner against BoomerAMG for 2D and 3D Poisson problems and showed how, in comparison to the 2D FMM, the more complicated oct-tree traversal for calculations of M2L and P2P kernels slows down the time-to-solution in the 3D FMM. Figure~\ref{fig:fmm_kernels} puts the 3D Helmholtz FMM in perspective with the 2D and 3D Poisson kernels. For the same problem size, the 3D Helmholtz is about an order of magnitude slower than the 3D Poisson because the Helmholtz operations are much more complicated to compute. Nevertheless, our model problems in~\sect \ref{sec:results} show that the FMM preconditioner requires small number of iterations to converge in comparison to the GMG, AMG, and IC preconditioners which may lead to a lower time-to-solution.\\

The model problem we consider in this section is defined on the unit cube $[0,1]^3$ and is characterized by homogeneous Dirichlet boundary conditions as follows:
\begin{subequations}
\begin{alignat}{3}
\nabla^2 u + \kappa^2 u & = 1 && \text{ in } & 
\domain,\label{eq:helm3_1}\\
u & = 0 && \text{ on } & \bound.
\end{alignat}
\label{eq:helm3}
\end{subequations}
with $\kappa = 0$.

\begin{figure}[t]
\centering
\includegraphics[width=0.7\textwidth,height=7cm]{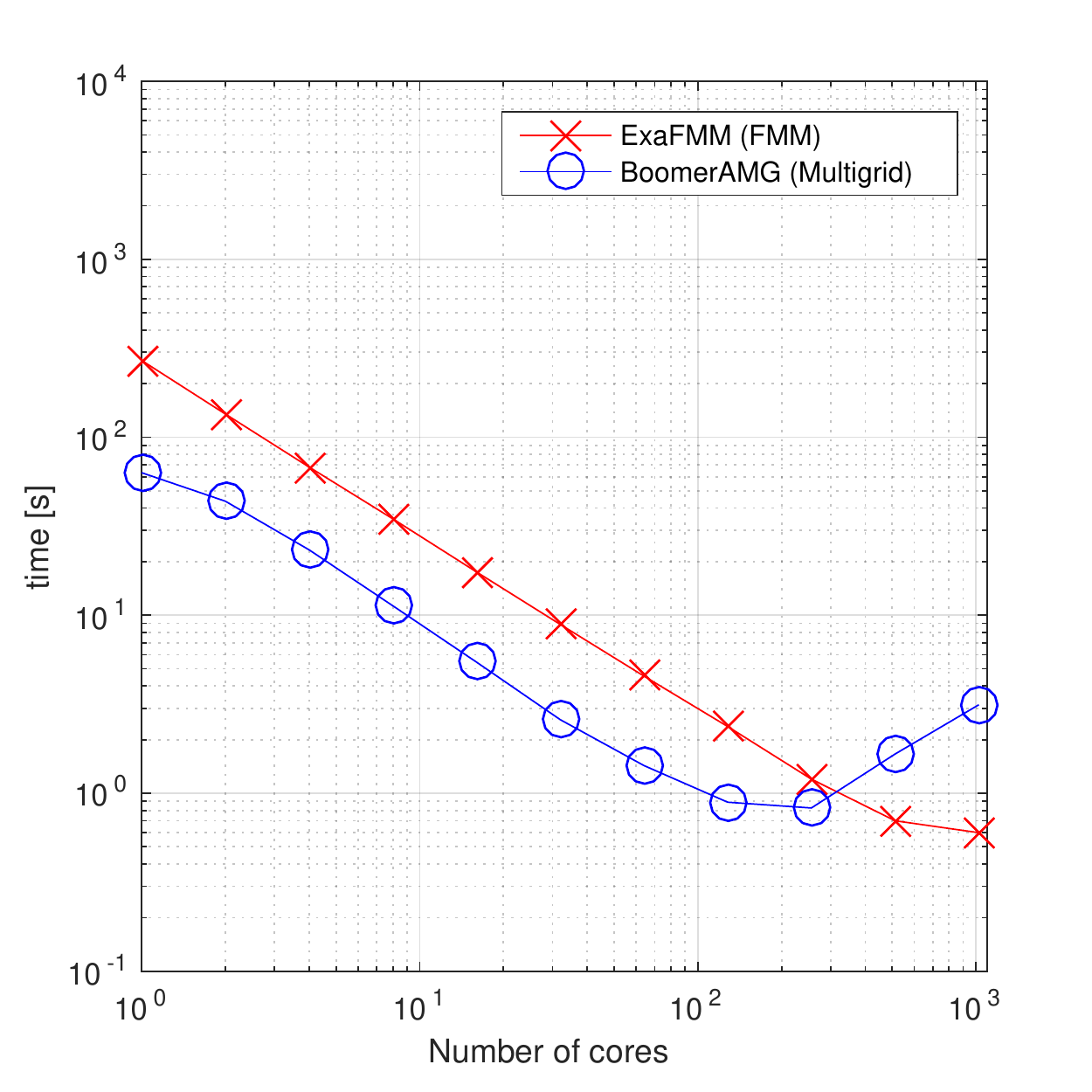}
\caption{Strong scaling of the FMM preconditioner.}
\label{fig:strong_scaling}
\end{figure}

Figure~\ref{fig:strong_scaling} shows the strong scalability of the FMM preconditioner on up to 1024 cores. The order of expansion for the FMM is set to $p=6$ and $\theta=0.4$, which gives about six significant digits of accuracy. The mesh size used in these runs is $N=4096^2$. Timings were measured with the PETSc run time option \texttt{-log\_summary}. We can see from Figure~\ref{fig:strong_scaling} that the FMM-based preconditioner strong scales quite well with respect to the problem size and number of cores.

\section{Conclusions}
\label{sec:conc}

Sparse linear solvers asymptotically dominate cost in scaling implicit mesh-based discretizations of PDEs. In Krylov-type solvers, preconditioning is usually the highest cost in terms of compute time, memory bandwidth, and memory size. Therefore, improving the scalability of preconditioning attacks the leading bottleneck in the scalability of mesh-based PDEs. The unique combination of $\mathcal{O}(N \log N)$ complexity and compute bound kernels makes the FMM an interesting candidate for preconditioners on future architectures with low Byte/flop ratios. In the present work, the FMM is employed as a preconditioner for Krylov subspace methods applied to discretizations of the Helmholtz equation prevalent in computational simulations. We tested the FMM-based preconditioner for the GMRES method and showed that the fast multipole method can be successfully coupled to the boundary element method to give an effective preconditioner (proper attention to the truncation error of the BEM relative to that of the PDE itself is needed). Our results show that the FMM-based preconditioner achieves both mesh-independent and wavenumber-independent convergence rate, for the tested values of $\kappa$ and $h$, while excelling in scalings on commodity architecture supercomputers. Compared with other methods exploiting the low rank structure of off-diagonal blocks, FMM-preconditioned Krylov iteration may reduce the amount of communication because it is matrix-free and exploits the tree structure of FMM.

\section{Acknowledgements}
The authors would like to acknowledge the open source software packages that made this work possible: PETSc~\cite{petsc-user-ref,petsc-web-page}, PetIGA~\cite{PetIGA} and IFISS~\cite{elman2007,ifiss}. We thank Jennifer Pestana and Stefano Zampini for helpful discussions and comments and Lisandro Dalcin for his help with the PetIGA framework. For computer time, this research used the resources of the Supercomputing Laboratory at King Abdullah University of Science \& Technology (KAUST) in Thuwal, Saudi Arabia.


\section*{References}
\bibliography{refs}

\end{document}